\documentclass[a4paper,11pt]{article}
\pdfoutput=1 % if your are submitting a pdflatex (i.e. if you have
\usepackage{jheppub} % for details on the use of the package, please
\usepackage[T1]{fontenc} % if needed

\usepackage{siunitx} % Required for alignment

\makeatletter
\def\@fpheader{\relax}
\makeatother

\title{Quantum Approach to Dehn Surgery Problem }

\author[]{Dongmin Gang}

\affiliation[]{Center for Theoretical Physics, Seoul National University}

\emailAdd{arima275@snu.ac.kr}

\abstract{We propose a new algorithm for Dehn surgery problem, finding exceptional Dehn filling slopes for a given hyperbolic 3-manifold with a torus boundary, using a quantum invariant called "3D index". The invariant is defined using an ideal triangulation of the cusped 3-manifold. We test  the algorithm for many examples. }

\begin{document} 
\maketitle
\flushbottom

\section{Introduction}
\label{sec:intro}

Quantum physics sometimes provide an efficient way of solving  mathematical problems in an unexpected way.  One famous example is the Shor's algorithm to the problem of integer factorization. Quantum knot invariants \cite{jones1985,witten1989quantum} are another good examples. They appearas basic physical observables in topological quantum field theories (TQFTs) and give  partial answer to the problem of distinguishing knots. 

In this note, we propose an efficient way to determine a basic property (hyperbolic or not) of a closed 3-manifold from its Dehn surgery representation. According to Lickorish-Wallace theorem \cite{Lickorish,Wallace}, every closed manifolds are given by a Dehn surgery along a link $L$ in $S^3$. Thurston proves that all knots but torus/satellite knots are hyperbolic \cite{thurston1979geometry}. He further prove that for given a hyperbolic 3-manifold $N$ with a torus boundary, there are only finite many exceptional slopes. A boundary slope, a primitive cycle $p\mu+q \lambda \in H_1 (\partial N ,\mathbb{Z})$ up to sign, is called `exceptional' if the closed manifold $M=N_{p \mu+ q \lambda}$ obtained by performing a Dehn filling along the slope is non-hyperbolic. Finding all the exceptional slopes for given $N$ is a non-trivial question and we call the problem as `Dehn surgery problem'. Classical approach to the problem is using an ideal triangulation of $N$. The hyperbolic structure (if exists) on the Dehn filled manifold can be obtained by gluing hyperbolic structures of each ideal tetrahedron in a compatible way. 
In this way, the topological problem can be translated into a problem of solving algebraic  equations, gluing equations, given   in \eqref{gluing eqns-1} and \eqref{gluing eqns-2} \cite{neumann1985volumes}.

In our approach, on the other hand,  we "quantize" the gluing equations  and find a wave-function satisfying  the quantum operator equations. A consistent quantization has been studied in \cite{Dimofte:2011gm,Bae:2016jpi} and the resulting wave functions turn out to be  powerful invariants. Different representations (labelled by an intger $k$) of a quantum torus algebra   give different wave-functions \cite{Dimofte:2014zga}. Refer to \cite{Dimofte:2012qj,Bae:2016jpi} for the wave-function with $k=1$ and its application to the volume conjecture \cite{1996q.alg.....1025K,1999math......5075M,Chen:2015wfa,Dimofte:2012qj,Gang:2017cwq}. The  wave-function associated to $k=0$  is called 3D index \cite{Dimofte:2011gm,Gang:2018wek} and turns out to be directly relevant to the Dehn surgery problem as we will see. 
The wave-function can be constructed by taking basic  operations, superposition and projection,  of quantum mechanics on the product of  tetrahedron's wave-functions $\mathcal{I}_{\Delta}$. See eq.~\eqref{3D index for N} and \eqref{Quantum Dehn filling} for the explicit expression. 
The final wave-function $\mathcal{I}_M(x)$  for Dehn filled closed manifold $M$ is generically an infinite power series in a   quantum parameter  $x^{1/2}$. 
 From a basic property of the power series, we can determine the  geometric structure of $M$ as proposed in  \eqref{Exceptional-index}. In the process, unlike the classical approach,  we do not need to solve  algebraic equations.

\section{Solving Dehn surgery problem using 3D index}
After a  summary of the recently proposed quantum 3-manifold invariant, 3D index \cite{Dimofte:2011py,Gang:2018wek}\footnote{3D index is firstly introduced in \cite{Dimofte:2011py} for 3-manifolds with torus boundaries and extended for closed 3-manifolds in \cite{Gang:2018wek} by incorporating Dehn filling operation. See also \cite{Garoufalidis:2016ckn} for an interesting relation between the 3D index and counting of normal surfaces.  },  we will propose a way of determining whether the 3-manifold is hyperbolic or not from a simple property of the invariant. We also conjecture that the invariant vanishes when the 3-manifold is Lens space.

\subsection{3D index : Review}
Let $N$ be a 3-manifold with a torus boundary. For given choice of an ideal triangulation $\mathcal{T}$ of the $N$, composed of $K$ ideal tetrahedrons with a semi-angle structure and a choice of basis  $(\mu,\lambda)$  of $H_1 (\partial N, \mathbb{Z})$, Thurston's gluing equations are  given in the following form\footnote{Refer to, for example, section 2 of \cite{Dimofte:2012qj} for details on the gluing equations.}
\begin{align}
\begin{split}
& z_i (z'_i) (z''_i)=-1\quad \textrm{and} \quad z_i^{-1}+z''_i-1=0\;,\; \quad \textrm{for } i=1,\ldots, K\;,
\\
&\prod_{j=1}^{K} z_j^{F_{ij}} (z'_j)^{F'_{ij}} (z''_j)^{F''_{ij}} = (-1)^2\;,\quad \textrm{for } i=2,\ldots, K\;,
\\
&\prod_{j=1}^{K} z_j^{M_{j}} (z'_j)^{M'_{j}} (z''_j)^{M''_{j}} = \mathbf{m}^2\;,
\\
&\prod_{j=1}^{K} z_j^{L_{j}} (z'_j)^{L'_{j}} (z''_j)^{L''_{j}}  = \ell^{2 }\;\;. \label{Guling matrices}
\end{split}
\end{align}
We choose a basis $(\mu, \lambda)$ such that $\lambda \in \textrm{Ker} [H_1 (\partial N,\mathbb{Z})\rightarrow H_1 (N, \mathbb{Z}_2)]$.
Here 
\begin{align}
F_{ij}, F'_{ij}, F''_{ij} \in \{ 0,1,2\}\;,  \quad M_j , M_j', M_j'', L_j, L_j', L_j'' \in \mathbb{Z}\;.
\end{align}
Integrating $(z_i',z_i'')$ using the equations in the first line, the equations can be simplified as
\begin{align}
\begin{split}
&\prod_{j=1}^{K} z_j^{A_{ij}} (1-z^{-1}_j)^{B_{ij}} = \mathbf{m}^{2 \delta_{i,1}}(-1)^{\nu_i}\;, \quad \quad i=1,\ldots, K\;,
\\
&\prod_{j=1}^{K} z_j^{2C_{j}} (1-z^{-1}_j)^{2D_{j}} = \ell^{2 }(-1)^{2 \nu_{K+1}}\;\;, \label{gluing eqns-1}
\end{split}
\end{align}
where
\begin{align}
\begin{split}
&A_{1j} = M_j - M_j'\;, \;  B_{1j} = M_j''-M_j'\;, \; \nu_1 = -\sum M_j'\;,
\\
&A_{ij} = F_{ij}- F'_{ij}\;, \; B_{ij} = F''_{ij}- F_{ij}'\;, \; \nu_i = 2- \sum F'_{ij}\;, \quad \textrm{for $i=2,\ldots, K$}\;,
\\
&2C_{j} = L_j- L'_j\;, \; 2D_j =L_j''-L_j\;, \; 2\nu_{K+1} = - \sum L'_j\;.
\end{split}
\end{align}
$\delta_{i,j}$ is the Kronecker delta symbol
\begin{align}
\delta_{i,j} = \begin{cases}
1\;,\; \textrm{ if $i=j$}\\
0\;,\; \textrm{ otherwise}
\end{cases}
\end{align}
A solution (if exists) to the gluing equations  with following additional conditions 
\begin{align}
\mathbf{m}^{2p} \ell^{2q}=1\;, \quad 0<\textrm{Im}[\log(z_i)]<\pi \label{gluing eqns-2}
\end{align}
could give a hyperbolic structure on $N_{p\mu +q \lambda}$. 
More generally,  solutions to the equations give $PGL(2,\mathbb{C})$ flat connections $\rho \in \textrm{Hom}\big{[}\pi_1 (N)\rightarrow PGL(2,\mathbb{C})\big{]}/(\textrm{conj})$ on $N$ with following boundary holonomies
\begin{align}
\begin{split}
\rho(\mu) = \begin{pmatrix} 
\mathbf{m} & * \\
0 & \mathbf{m}^{-1} 
\end{pmatrix}
\;, \quad
\rho(\lambda) = \begin{pmatrix} 
\ell & * \\
0 & \ell^{-1} 
\end{pmatrix}
\end{split}
\end{align}
Due to a symplectic property of the gluing equations \cite{neumann1985volumes}, the small matrices can be completed to a $Sp(2K,\mathbb{Z})$ matrix $G$ whose first $(K+1)$-rows are given as follows
\begin{align}
(G)_{i,j}= A_{ij}\;, \; (G)_{i,K+j}= B_{ij}\;, \;(G)_{K+1,i} = C_i\;, \; (G)_{K+1,K+i}=D_i\;.
\end{align}
There is a choice of $(K-1)$-rows of $G$ such that  the resulting matrix $G$ becomes a symplectic matrix. Using a particular choice of $G$, the 3D index $\mathcal{I}^{(\mu, \lambda)}_N$ for a 3-manifold $N$ with a torus boundary  is defined as \cite{Dimofte:2011py}
\begin{align}
\begin{split}
&\textrm{{\bf Def} : 3D index for $N$ with resect to  a basis choice $(\mu,\lambda)\in H_1 (\partial N, \mathbb{Z})$ }
\\
&\mathcal{I}_{N}^{(\mu,\lambda)}\;:\;  (m,e) \in \mathbb{Z}\times \mathbb{Z} \quad  \rightarrow  \quad  \mathcal{I}_{N}^{(\mu,\lambda)} (m,e;x) \in \mathbb{Z}((x^{1/2}))\;,
\\
&\textrm{where}
\\
&\mathcal{I}_{N}^{(\mu,\lambda)} (m,e;x) 
\\
&=\sum_{(e_2,\ldots,e_{K}) \in \mathbb{Z}^{K-1}} \left( (-x^{\frac{1}2})^{\langle \nu, \gamma \rangle } \prod_{i=1}^K  \mathcal{I}_{\Delta} \big{(}(G^{-1} \cdot \gamma)_i , (G^{-1} \cdot \gamma)_{K+i} ;x \big{)} \right)\bigg{|}_{m_1 \rightarrow m, \;e_1\rightarrow e, \;m_{I>1} \rightarrow 0} \;, 
\\
&\textrm{where } \gamma := (m_1,\ldots, m_K, e_1, \ldots, e_K)^t\; \textrm{and}\; \langle \nu, \gamma \rangle :=  \sum_{i=1}^K\nu_i e_i - \nu_{K+1}m_1 \;. \label{3D index for N}
\end{split}
\end{align}
The expression is independent on the choice of the last $(K-1)$-rows of $G$ as far as the matrix is symplectic. Here, the tetrahedron index  $\mathcal{I}_{\Delta}(m,e;x)$ is given by 
\begin{align}
\begin{split}
&\mathcal{I}_{\Delta} (m,e;x) = \sum_{n=[e]}^\infty \frac{(-1)^n x^{\frac{1}2 n (n+1)-(n+\frac{1}2 e)m}}{(x)_n(x)_{n+e}} \;,
\\
&\textrm{where $[e]:=\frac{1}2 (|e|-e)$ and $(x)_n:= (1-x)(1-x^2)\ldots (1-x^n)$}.
\label{tetrahdron index}
\end{split}
\end{align}
The tetrahedron index can be consider as a formal infinite power series in $x^{1/2}$ by interpretating
\begin{align}
\frac{1}{1-x^n} = \sum_{k=0}^{\infty} x^{nk}=1+x^n+x^{2n}+\ldots\;.
\end{align}
For example,
\begin{align}
\begin{split}
&\mathcal{I}_{\Delta}(0,0;x) = 1-x-2x^2-2x^3-2x^4+x^6+5 x^7+\ldots\;,
\\
&\mathcal{I}_{\Delta}(1,0;x) = -x-x^2+x^4+3x^5+4x^6+6 x^8+6 x^9+\ldots\;,
\\
&\mathcal{I}_{\Delta}(0,1;x) = 1-x^2-2x^3-3x^4-3x^5-x^7+x^8+5 x^9+ \ldots\;,
\\
&\mathcal{I}_{\Delta}(1,1;x) = -x^{\frac{3}2}-x^{\frac{5}2}-x^{\frac{7}2}+x^{\frac{11}2}+3 x^{\frac{13}2}+5 x^{\frac{15}2}+7 x^{\frac{17}2}+\ldots\;.
\end{split}
\end{align}
The infinite summation in \eqref{3D index for N} is proven to converge\footnote{The convergence means that if you want to compute the index up to $o(x^{\sharp})$ with a fixed $\sharp \in \mathbb{Z}/2$ you only need to sum over finite many $(e_2,\ldots, e_K)$ and other contributions will only give higher power in $x$.  } when the triangulation $\mathcal{T}$ has a semi-angle structure \cite{garoufalidis20151}. 
The index is proven to be invariant under local 2-3 move of $\mathcal{T}$ if the semi-angle structure are preserved under the move. Although there is no proof that the index is independent on the choice of $\mathcal{T}$ with seim-angle structure, we will assume it is true. The index for general $N$ has following $\mathbb{Z}_2$-symmetry 
\begin{align}
\mathcal{I}_N^{(\mu, \lambda)} (m,e;x)=\mathcal{I}_N^{(\mu, \lambda)} (-m,-e;x)\;.
\end{align}

\paragraph{Example : $N=S^3\backslash \mathbf{4_1}$} According to SnapPy \cite{SnapPy}, the knot complement can be triangulated by two tetrahedra. The gluing matrices \eqref{Guling matrices} are
\begin{align}
\begin{split}
&(M_1, M_1',M_1'',M_2 ,M_2',M_2'') = (1,0,0,0,-1,0)\;,
\\
&(L_1, L_1',L_1'',L_2 ,L_2',L_2'') = (0,0,0,0,-2,2)\;,
\\
&(F_{21}, F_{21}',F_{21}'',F_{22} ,F_{22}',F_{22}'') = (0,1,2,1,2,0)\;,
\\
&\Rightarrow
\\
& (A)_{ij} = \begin{pmatrix} 
1 & 1 \\
-1 & -1 
\end{pmatrix} \;, \quad (B)_{ij} = \begin{pmatrix} 
0 & 1 \\
1 & -2 
\end{pmatrix}\;,\; 
\\
&(C_1,C_2 )= (0,1)\;, \; (D_1, D_2)= (0,2) \;, \; (\nu_1, \nu_2, \nu_3) = (1,-1,1)\;.
\end{split}
\end{align}
The matrices $(A,B,C,D)$ can be completed to a following $Sp(4,\mathbb{Z})$ matrix
\begin{align}
G = \begin{pmatrix} 
1 & 1 & 0 & 1\\
-1 & -1 & 1 & -2 \\
0 & 1 & 0 & 2 \\
-1 & 0 & 0 & 0\\
\end{pmatrix}
\end{align}
Then, using \eqref{3D index for N}, the index is given by\footnote{When $N$ is a  knot complement on $S^3$, there is canonical basis of $H_1 (\partial N, \mathbb{Z})$, meridian $\mu$ and longitude $\lambda$. }
\begin{align}
&\mathcal{I}_{S^3\backslash \mathbf{4_1}}^{(\mu, \lambda)} (m,e;x) =  \sum_{e_2 \in \mathbb{Z}}(-q^{\frac{1}2})^{-m-e_2}\mathcal{I}_\Delta (-e_2, -e_2;x) \mathcal{I}_{\Delta}(2e_2 +2m,-e_2-m;x)\;. \label{Ind 41}
\end{align}
Listing some,
\begin{align}
\begin{split}
&\mathcal{I}_{S^3\backslash \mathbf{4_1}}^{(\mu, \lambda)} (0,0;x) = 1-2x-3x^2+2x^3+ 8 x^4+18 x^6+18 x^6+14 x^7+\ldots\;,
\\
&\mathcal{I}_{S^3\backslash \mathbf{4_1}}^{(\mu, \lambda)} (0,\pm 1;x) = -2x-2x^2+2x^3+8 x^4+ 16 x^6+ 16 x^6+10 x^7+\ldots\;,
\\
&\mathcal{I}_{S^3\backslash \mathbf{4_1}}^{(\mu, \lambda)} (\pm 1,0;x) = -2x^{\frac{3}2} + 4 x^{\frac{7}2}+10 x^{\frac{9}2}+ 14 x^{\frac{11}2} + 10 x^{\frac{13}2}-2 x^{\frac{15}2}+\ldots\;,
\\
&\mathcal{I}_{S^3\backslash \mathbf{4_1}}^{(\mu, \lambda)} (\pm 1,\pm 1;x) = -x-x^2 +2 x^3+7 x^4+11 x^5 +11 x^6 + 3 x^7+\ldots\;.
\end{split}
\end{align}
\\
\paragraph{Example : $N=S^3\backslash \mathbf{5_2}$} According to SnapPy, the knot complement can be triangulated by three tetrahedra. Using the ideal triangulation, the index is given by 
\begin{align}
\begin{split}
&\mathcal{I}_{S^3\backslash \mathbf{5_2}}^{(\mu, \lambda)} (m,e;x) =  \sum_{e_2 ,e_3 \in \mathbb{Z}}(-x^{\frac{1}2})^{-(e_3+m)} \mathcal{I}_{\Delta}(e_3,e_2;x) \mathcal{I}_\Delta(e_3+m,-e-2 e_2-e_3-2 m;x) 
\\
&\qquad \qquad \qquad \qquad \qquad  \times \mathcal{I}_\Delta(e_3+2 m,e+e_2+m;x)\;. \label{Ind 52}
\end{split}
\end{align}
Listing some,
\begin{align}
\begin{split}
&\mathcal{I}_{S^3\backslash \mathbf{5_2}}^{(\mu, \lambda)} (0,0;x) = 1-4x-x^2+16x^3+ 26 x^4+ 23 x^5-34 x^6 -122  x^7+\ldots\;,
\\
&\mathcal{I}_{S^3\backslash \mathbf{5_2}}^{(\mu, \lambda)} (0,\pm 1;x) = -3x+14x^3+ 22 x^4+ 16 x^5- 36 x^6-116 x^7+\ldots\;,
\\
&\mathcal{I}_{S^3\backslash \mathbf{5_2}}^{(\mu, \lambda)} (\pm 1,0;x) = -x^{\frac{1}2} + 5 x^{\frac{5}2}+12 x^{\frac{7}2}+ 10 x^{\frac{9}2} - 15 x^{\frac{11}2}-64 x^{\frac{13}2}+\ldots\;,
\\
&\mathcal{I}_{S^3\backslash \mathbf{5_2}}^{(\mu, \lambda)} ( 1, 1;x) =\mathcal{I}_{S^3\backslash \mathbf{5_2}}^{(\mu, \lambda)} (-1, -1;x)= 2x^2+ 4 x^3+ 5x^4- 7 x^5-31 x^6-69 x^7-95 x^8+\ldots\;,
\\
&\mathcal{I}_{S^3\backslash \mathbf{5_2}}^{(\mu, \lambda)} ( 1, -1;x) =\mathcal{I}_{S^3\backslash \mathbf{5_2}}^{(\mu, \lambda)} (-1, 1;x) = -2x+x^2+12 x^3+17 x^4+8 x^5-38 x^6 -107 x^7-186 x^8+\ldots\;.
\end{split}
\end{align}
\\
The 3D index for a closed 3-manifold $N_{p\mu+q \lambda}$, obtained by performing Dehn filling on $N$ along a boundary cycle $p\mu+q\lambda$, is defined as follows (Appendix A of \cite{Gang:2018wek}):
\begin{align}
\begin{split}
&\textrm{{\bf Def} : 3D index for $N_{p\mu + q\lambda}$ }
\\
& \mathcal{I}_{N_{p\mu+q\lambda}}(x)  \in \mathbb{Z}((x^{1/2})) \textrm{ is defined by follwoing infinte sum}
\\
&\mathcal{I}_{N_{p\mu+q\lambda}}(x) =\sum_{(m,e)\in \mathbb{Z}^2} \mathcal{K}(m,e;p,q;x)\mathcal{I}^{(\mu,\lambda)}_N(m,e;x)\;, \quad \textrm{where}
\\
&\mathcal{K}(m,e;p,q;x) =\frac{1}2 (-1)^{rm+2se} \bigg{(}\delta_{p m+2q e,0} (x^{\frac{r m+2s e}2}+ x^{-\frac{r m+2s e}2}) - \delta_{p m+2qe,-2}-\delta_{p m+2q e,2} \bigg{)}\;. \label{Quantum Dehn filling}  
\end{split}
\end{align}
Here $(p,q)$ are coprime and $(r,s)$ are integers satisfying $qr-ps=1$. The choice of $(r,s)$ is not unique but can be shifted by $(r,s)\rightarrow (r+p, s+q)$.  The above expression is invariant under the shift.  To show the topological invariance of the index, we need to show that
\begin{align}
\textrm{If $N_{p\mu+ q \nu}$ is homeomorphic to $\tilde{N}_{\tilde{p}\tilde{\mu}+ \tilde{q} \tilde{\nu}}$, then $\mathcal{I}_{N_{p\mu+q\lambda}}(x)=\mathcal{I}_{\tilde{N}_{\tilde{p}\tilde{\mu}+\tilde{q}\tilde{\lambda}}}(x)$}.
\end{align}
Again, we currently do not have a proof of the invariance. The above formula in eq.~\eqref{Quantum Dehn filling} is derived by applying a general  method in TQFTs   to a particular case and the invariance is very likely. For many examples, the invariance has been checked \cite{Gang:2018wek}. 
\paragraph{Example :} From a surgery calculus, one can check that
\begin{align}
(S^3\backslash \mathbf{4_1})_{-5\mu+\lambda} = (S^3\backslash \mathbf{5_2})_{5\mu+\lambda}\;.
\end{align}
It is compatible with following computation (using eq.~\eqref{Ind 41} and \eqref{Ind 52})
\begin{align}
\begin{split}
&\sum_{(m,e)\in \mathbb{Z}^2} \mathcal{K}(m,e;5,1;x) \mathcal{I}^{(\mu,\lambda)}_{S^3 \backslash \mathbf{5_2}}(m,e;x) =  \sum_{(m,e)\in \mathbb{Z}^2} \mathcal{K}(m,e;-5,1;x) \mathcal{I}^{(\mu,\lambda)}_{S^3 \backslash \mathbf{4_1}}(m,e;x) 
\\
&=1-x-2x^2-x^3-x^4+x^5+2 x^6+7 x^7+8 x^8+ 12 x^9+\ldots
\end{split}
\end{align}
\\

\subsection{Detecting exceptional slopes using 3D Index}
We propose following property of the 3D index 
\begin{align}
\begin{split}
&\textrm{\bf Conjecture }
\\
&\mathcal{I}_{M=N_{p\mu+q\lambda}}(x) = 
\begin{cases}
\textrm{Infinite series starting with 1+\ldots }\;, \; \textrm{if $M$ is hyperbolic}
\\
0,1\textrm{ or }  \infty \textrm{ (does not converge)}\;, \;\textrm{if $M$ is non-hyperbolic}\;\\
\end{cases}
\\
&\textrm{If $M$ is a Lens space, then $I_{M}(x)=0$}\;.
\label{Exceptional-index}
\end{split}
\end{align}
There is a partial understanding of the conjecture from the physics of M5-branes  wrapped on 3-manifolds \cite{Gang:2018wek}.   
Using the proposed property we can determine exceptional slopes by computing the index in the power series expansion of $x^{1/2}$. 
\\
\\
Now let us give examples supporting the conjectures.
\\
\paragraph{Example : Exceptional slopes for $S^3\backslash \mathbf{4_1}$} Using the formulae in eq.~\eqref{Ind 41} and \eqref{Quantum Dehn filling}, we can compute the index $I_{M= (S^3\backslash \mathbf{4_1})_{p\mu + q\lambda}}(x)$ in series of $x^{\frac{1}2}$. The 3D indices for $|p|+|q|\leq 7$ are \footnote{"0" ("1") in the table means that we checked that the index is  $0+o(x^{\sharp})$ ($1+o(x^{\sharp})$) up to  sufficiently high  order $\sharp$. From experiment, we confirm that 1) if the coefficient of $x^0$  vanishes, then coefficients of higher orders also vanish and  2) if the index is $1+o(x^{\sharp})$ up to sufficiently high $\sharp$ (for every cases we checked, $\sharp =3$ is sufficient), then the index remains as $1$ up to even higher order.}
\begin{align}
\begin{tabular}{ l| c }
  $(p,q)$ & $\mathcal{I}_{(S^3\backslash \mathbf{4_1})_{p\mu+q \lambda}}(x)$\\ \hline
  $(1,0)$ & $0$  \\
  $(0,1)$ & $1$  \\
  $(1,\pm 1)$ & $1$  \\
  $(1,\pm 2)$ & $1-2x-3x^2+3 x^4+10 x^5+\ldots$  \\
  $(2,\pm 1)$ & $1$  \\
  $(1,\pm 3)$ & $1-2x-3x^2 + x^3+6x^4+13 x^5+\ldots$  \\
  $(3,\pm 1)$ & $1$  \\
  $(1,\pm 4)$ & $1-2x-3x^2+x^3+6 x^4+13 x^5+\ldots$  \\
  $(2,\pm 3)$ & $1-2x-2x^{\frac{3}2}-3x^2+x^3+4 x^{\frac{7}2}+\ldots$  \\
  $(3,\pm 2)$ & $1-2x-3x^2-x^3+3 x^4+13 x^5+\ldots$  \\  
  $(4,\pm 1)$ & $\infty$  \\
  $(1,\pm 5)$ & $1 - 2 x- 3 x^2+  x^3+6 x^4 +\ldots$  \\
  $(5,\pm 1)$ & $1 -  x- 2 x^2- x^3- x^4 +x^5+\ldots$  \\
  $(1,\pm 6)$ & $1 -  2x- 3 x^2+ x^3+ 6 x^4 + 13 x^5+\ldots$  \\
  $(2,\pm 5)$ & $1 -  2x-2x^{\frac{3}2}- 3 x^2+2 x^3+ 4 x^{\frac{7}2} +\ldots$  \\
  $(3,\pm 4)$ & $1 -  2x- 3 x^2+ x^3+ 5 x^{4} +15 x^5+\ldots$ \\
  $(4,\pm 3)$ & $1 -  3x- 4 x^2+ 3x^3+ 12 x^{4} +27 x^5+\ldots$ \\
  $(5,\pm 2)$ & $1 -  2x- 4 x^2- x^3+ 5 x^{4} +17 x^5+\ldots$ \\
  $(6,\pm 1)$ & $1 -  x^{\frac{1}2}-  x^{\frac{3}2}- x^2- x^{\frac{5}2} + x^{\frac{9}2}+  3 x^5 +\ldots$
  \\
  \hline
\end{tabular}
\end{align}
By applying the criterion in eq.~\eqref{Exceptional-index}, we find following 10 exception slopes 
\begin{align}
\{\textrm{Exceptional slopes of $S^3\backslash \mathbf{4_1}$}\} = \{1/0, 0, \pm 1, \pm 2, \pm 3, \pm 4 \}\;.
\end{align}
This is compatible with known result in literature, for example \cite{gordon1998dehn}.  Further, applying the criterion in eq.~\eqref{Exceptional-index}, we have following candidate for Lens space slopes
\begin{align}
\{\textrm{Candidates for Lens space slopes of $S^3\backslash \mathbf{4_1}$}\} = \{1/0 \}\;.
\end{align}
Obviously,  $(S^3\backslash \mathbf{4_1})_{\mu} = S^3$ is a Lens space. 
\\
\paragraph{Example : Exceptional slopes for $S^3\backslash \mathbf{5_2}$} Using the formulae in eq.~\eqref{Ind 52} and \eqref{Quantum Dehn filling}, we can compute the index $I_{M= (S^3\backslash \mathbf{5_2})_{p\mu + q\lambda}}(x)$ in series of $x$. Listing the 3D indices for $|p|+|q|\leq 6$
\begin{align}
\begin{tabular}{ l| c }
  $(p,q)$ & $\mathcal{I}_{(S^3\backslash \mathbf{5_2})_{p\mu+q \lambda}}(x)$\\ \hline
  $(1,0)$ & $0$  \\
  $(0,1)$ & $\infty $ \\
  $(1,- 1)$ & $1-2 x-3x^2+3 x^4+ 10 x^5+\ldots$  \\
  $(1,1)$ & $1$  \\
  $(1,- 2)$ & $1-5x-x^2+18 x^3+30 x^4+\ldots$  \\
  $(1, 2)$ & $1-5x-3x^2+16 x^3+32 x^4+\ldots$  \\
  $(2,- 1)$ & $1-x^{\frac{1}2} - x -2 x^{\frac{3}2} - x^2 + 3x^3+\ldots$  \\
  $(2,1)$ & $1$  \\
  $(1,- 3)$ & $1-5x-x^2 +16 x^3+ 30 x^4+34 x^5+\ldots$  \\
  $(1, 3)$ & $1-5x-x^2 +16 x^3+ 28 x^4+32 x^5+\ldots$  \\
  $(3,- 1)$ & $1-x-x^2+2x^3+5 x^4+9 x^5+\ldots$  \\
  $(3, 1)$ & $1$  \\
  $(1, -4)$ & $1-5x-x^2+16x^3+30 x^4+32 x^5+\ldots$  \\
  $(1, 4)$ & $1-5x-x^2+16x^3+30 x^4+32 x^5+\ldots$  \\
  $(2,- 3)$ & $1-x^{\frac{1}2}-4x-x^2+5x^{\frac{5}2}+18 x^3+\ldots$  \\
  $(2,3)$ & $1-x^{\frac{1}2}-4x-3x^2+7x^{\frac{5}2}+14 x^3+\ldots$  \\
  $(3,- 2)$ & $1-4x-x^2+18x^3+29 x^4+28 x^5+\ldots$  \\ 
  $(3, 2)$ & $1-4x-4x^2+9x^3+25 x^4+45 x^5+\ldots$  \\ 
  $(4,- 1)$ & $1 - x+4 x^3+ 11 x^4+13 x^5+\ldots$  \\
  $(4, 1)$ & $\infty$  \\
  $(1,- 5)$ & $1 - 5 x-  x^2+ 16 x^3+30 x^4 +\ldots$  \\
  $(1,5)$ & $1 - 5 x-  x^2+ 16 x^3+30 x^4 +\ldots$  \\
  $(5,-1)$ & $1 -  x-  x^2+ 2 x^3+4 x^4 + 7 x^5+\ldots$  \\
  $(5,1)$ & $1 -  x- 2 x^2- x^3- x^4 +  x^5+\ldots$  \\
  \hline
\end{tabular}
\end{align}
By applying the criterion in eq.~\eqref{Exceptional-index}, we find following 6 exceptional slopes which is compatible with known result \cite{gordon1998dehn}. 
\begin{align}
\{\textrm{Exceptional slopes of $S^3\backslash \mathbf{5_2}$}\} = \{1/0, 0,  1,  2,  3,  4 \}\;.
\end{align}
Further, applying the criterion in eq.~\eqref{Exceptional-index}, we have following candidate for Lens space slopes
\begin{align}
\{\textrm{Candidates for Lens space slopes of $S^3\backslash \mathbf{5_2}$}\} = \{1/0 \}\;.
\end{align}
Obviously,  $(S^3\backslash \mathbf{5_2})_{\mu} = S^3$ is a Lens space. 
\paragraph{Example : Exceptional slopes for $m003$} $N=m003$ denotes the 3-manifold called `sister of figure-eight knot complement'. The manifold is not a knot complement on $S^3$ and there is no canonical basis choice $(\mu, \lambda)$ of $H_1 (\partial N, \mathbb{Z})$. SnapPy  choose a particular basis of $H_1(\partial N, \mathbb{Z})$, which we denote $(A,B)$. Since the $A  \in \textrm{Ker} [H_1 (\partial N,\mathbb{Z})\rightarrow H_1 (N, \mathbb{Z}_2)]$, we choose
\begin{align}
\mu = B \;, \; \lambda = -A\;.
\end{align}
In the basis choice, the index for $m003$ is (using the ideal triangulation provided in SnapPy)
\begin{align}
\mathcal{I}_{m003}^{(\mu, \lambda)} = \sum_{e_2 \in \mathbb{Z}} (-x^{\frac{1}2})^{m-2e_2} \mathcal{I}_{\Delta}(e_2, m+e-2e_2;x) \mathcal{I}_{\Delta}(e_2-m,-e+m-2e_2;x)
\end{align}
Computing the index for $(m003)_{pA+q B}$ using \eqref{Quantum Dehn filling}, we find following 8 exceptional slopes
\begin{align}
\begin{tabular}{ l| c }
  $(p,q)$ & $\mathcal{I}_{(m003)_{pA+q B}}(x)=\mathcal{I}_{(m003)_{-p\lambda+q\mu}}(x)$\\ \hline
  $(0,1)$ & $0$  \\
  $(1,0)$ & $0$  \\
  $(1,-1)$ & $0$  \\
  $(1,-2)$ & $1$  \\
  $(1,1)$ & $1$  \\
  $(1,2)$ & $\infty$  \\
  $(2,-1)$ & $1$  \\
   $(3,-2)$ & $\infty$  \\
   \hline
  \end{tabular}
\end{align}
\begin{align}
\{\textrm{Candidates for Lens space slopes of $ m003$}\} = \{0,1/0,-1 \}\;.
\end{align}
Actually these all correspond to Lens spaces \cite{martelli2002dehn}. 

\acknowledgments
The contents of this paper was presented in QMAP seminar at UC Davis. We thank the audience for feedback.  We also thank Bruno Martelli and Ian Agol for kind email correspondence. The work is supported by Samsung
Science and Technology Foundation under Project Number
SSTBA1402-08.

%This is the most common positions for acknowledgments. A macro is
%available to maintain the same layout and spelling of the heading.


\bibliographystyle{JHEP}
\bibliography{ref}
\end{document}